\documentclass[a4paper,11pt]{article}

\usepackage[latin1]{inputenc}
\usepackage{amsmath}
\usepackage[francais]{babel}
\usepackage{amsfonts}
\usepackage{amssymb}
\usepackage{qsymbols}
\usepackage{euscript}
\usepackage{graphicx}

\newcounter{compt}
\begin{document}

\newtheorem{lem}{Lemma}
\newtheorem{thm}{Theorem}
\newtheorem{prop}{Proposition}
\newtheorem{Def}{Definition}


\newcommand{\N}{\mathbb{N}}
\newcommand{\R}{\mathbb{R}}
\newcommand{\Z}{\mathbb{Z}}
\newcommand{\C}{\mathbb{C}}
\newcommand{\K}{\mathbb{K}}
\newcommand{\Q}{\mathbb{Q}}
\newcommand{\M}{\mathbb{M}}
\newcommand{\A}{\mathbb{A}}
\newcommand{\Hm}{\mathbb{H}}
\newcommand{\Sy}{\mathbb{S}}
\newcommand{\mB}{\mathcal{B}}
\newcommand{\Prob}{\mathcal{P}}
\newcommand{\E}{\mathbb{E}}


\newcommand{\vp}{\varphi}
\newcommand{\ve}{\varepsilon}
\newcommand{\la}{\lambda}


\newcommand{\Part}[2]{\displaystyle\frac{\partial{#1}}{\partial{#2}}}
\newcommand{\Dd}[2]{\displaystyle\frac{\text{d}{#1}}{\text{d}{#2}}}
\newcommand{\Lun}{L^1(\R)}
\newcommand{\supp}{\text{supp}}
\newcommand{\cod}{\text{codim}}
\newcommand{\Ran}{\text{Ran}}
\newcommand{\Ker}{\text{Ker}}
\newcommand{\ch}{\text{ch}}
\newcommand{\sh}{\text{sh}}
\newcommand{\etoile}{\begin{center}$\qquad\qquad\ast\qquad\ast\qquad\ast\qquad\ast\qquad\ast\qquad\ast\qquad\ast$\end{center}}
\newcommand{\dt}{\text{d}}
\newcommand{\Filt}{\text{\cls F}\;}
\newcommand{\rayer}[1]{\setbox1=\hbox{#1}\mbox{}\rlap{\rule[.3ex]{\wd1}{1pt}}\unhbox1}
\newcommand{\der}{\partial}


\newcommand{\dis}{\displaystyle}


\newcommand{\indi}{1\!\!I}

\newcommand{\ind}{\begin{picture}(2,10)
1\put(-5.5,0.4){ \framebox(2.0,5.77){ }}
\end{picture}$\ $}
\newcommand{\inddeux}{\text{\begin{picture}(2,10)
1\put(-5.5,0.4){ \framebox(2.0,5.77){ }}
\end{picture}} }


\title{Global Existence of classical solutions\\ for a class of reaction-diffusion systems}
\author{El-Haj Laamri\\
Institut Elie Cartan\\
Universit\'e Henri Poincar\'{e}, Nancy 1\\
B.P. 239\\
54 506 Vandoeuvre-l\`es-Nancy\\
Elhaj.Laamri@iecn.u-nancy.fr}
\date{ Version du mardi 22 f\'{e}vrier \`{a}  9h}

\maketitle

{\bf Abstract } In this paper, we use duality arguments "\`{a} la  Michel Pierre" to establish  global existence of  classic solutions for a class of parabolic  reaction-diffusion systems modeling, for instance, the evolution of reversible chemical reactions.
\section{Introduction}
This paper is motivated by the general question of global existence in time of solutions to the  following reaction-diffusion system
$$(\mathcal S)\left\lbrace\begin{array}{lllll}
u_t-d_1\Delta u & = & w^\gamma-u^\alpha v^\beta & (0,+\infty)\times\Omega, &(E_1)\\
v_t-d_2\Delta v & = &  w^\gamma-u^\alpha v^\beta & (0,+\infty)\times\Omega, &(E_2)\\
w_t-d_3\Delta w & = & - w^\gamma+ u^\alpha v^\beta & (0,+\infty)\times\Omega, &(E_3)\\
\Part{u}{n}(t,x)=\Part{v}{n}(t,x)=\Part{w}{n}(t,x)    & = & 0 & (0,+\infty)\times\partial\Omega,\\
u(0,x) & = & u_0(x)\geq 0 & x\in\Omega,\\
v(0,x) & = & v_0(x)\geq 0 & x\in\Omega,\\
w(0,x) & = & w_0(x)\geq 0 & x\in\Omega,
\end{array}
\right.$$
 where  $\Omega$ is a bounded regular open subset of  $\mathbb{R}^N$, $(d_1, d_2, d_3,\alpha,\beta,\gamma  )\in (0,+\infty)^3\times[1,+\infty)^3$.
\vskip3mm
Note that the system $(\mathcal S)$ satisfies two main properties, namely  :\\
$(P)$ the nonnegativity of solutions of $(\mathcal S)$  is preserved for all time ;\\
$(M)$ the total mass of the components $u$, $v$, $w$  is a priori bounded on all
finite intervals $(0,t)$.
\vskip5mm
If $\alpha,\beta$ and  $\gamma$ are positive integers, system $(\mathcal S)$ is intended to describe  for example the evolution of a reversible chemical reaction of type
 $$\alpha U+\beta V\rightleftharpoons \gamma W$$
where $u$, $v$, $w$ stand for the density of $U$,   $V$ and $W$ respectively.\\
This chemical reaction is typical of general reversible reactions and contains
 the major difficulties encountered in a large class of similar problems as regards global existence of solutions.
\vskip3mm
Let us make precise what we mean by solution.
\vskip2mm
By {\it classical solution} to $(\mathcal S)$ on $Q_T=(0,T)\times \Omega$, we mean that, at least\\
$(i)$ $(u,v,w)\in \mathcal{C}([0,T);L^1(\Omega)^3)\cap L^\infty([0,\tau]\times\Omega)^3, \forall \tau \in (0,T)$ ;\\
$(ii)$ $\forall k,\ell=1\dots N$, $\forall p\in (1,+\infty)$
 $$\der_t u,\der_t v,\der_t w,\der_{x_k}u,\der_{x_k}v,\der_{x_k}w,\der_{x_kx_\ell}u,\der_{x_kx_\ell}v,\der_{x_kx_\ell}w, u,v,w \in L^p((0,T)\times\Omega) \,;$$
$(iii)$ equations in $(\cal S)$  are satisfied a.e (almost everywhere).
\vskip2mm
By {\it weak solution} to $(\mathcal S)$ on $Q_T=(0,T)\times \Omega$, we essentially  mean solution in the sense of distributions or, equivalently here, solution in the sens of the variation of constants formula with the corresponding semigroups. More precisely
\begin{eqnarray*}
u(t)&=&S_{d_1}(t)u_0 + \int_0 S_{d_1}(t-s)( w^\gamma(s)-u^\alpha(s) v^\beta(s))\,ds\\
v(t)&=&S_{d_2}(t)v_0 + \int_0 S_{d_2}(t-s)( w^\gamma(s)-u^\alpha(s) v^\beta(s))\,ds\\
w(t)&=&S_{d_3}(t)u_0 + \int_0 S_{d_3}(t-s)( -w^\gamma(s)+u^\alpha(s) v^\beta(s))\,ds
\end{eqnarray*}
where $S_{d_i}(.)$ is the semigroup generated in $L^1(\Omega)$  by $-d_i\Delta$  with homogeneous Neumann boundary condition, $1\leq i\leq 3$.
\vskip5mm
By just integrating the sum  $(E_1)+ (E_2)+ 2 (E_3)$ in space and time,  and taking into account the boundary conditions $\left(\dis\int_\Omega \Delta(d_1u+d_2v+d_3w)=0\right)$, we obtain
\begin{eqnarray}\label{Conservationofthetotalmass}
\int_\Omega u(t)+v(t)+2w(t) = \int_\Omega u_0+v_0+2w_0\quad\quad t\geq 0.
\end{eqnarray}
Together with the nonnegativity of $u$, $v$ and $w$, estimate (\ref{Conservationofthetotalmass}) implies that
\begin{eqnarray}\label{estimateLUN}
\forall t\geq0\;,\; \|u(t)\|_{L^1(\Omega)}, \|v(t)\|_{L^1(\Omega)}, \|w(t)\|_{L^1(\Omega)}\leq \|u_0+v_0+2w_0\|_{L^1(\Omega)}.
\end{eqnarray}
In other words, the total mass of three components does not blow up ;  $u(t)$, $v(t)$ and  $w(t)$ rest bounded in $L^1(\Omega)$  uniformly in time.
\vskip1mm
Although one has uniform $L¹$-bound in time, classical solutions may not globally exist for  diffusion coefficients $d_1$, $d_2$, $d_3$ which are not equal (global existence obviously holds if $d_1=d_2=d_3$). As surprisingly proved in \cite{PS} and \cite{S2}, it may indeed happen that, under assumptions $(P)$ and $(M)$, solutions blow up in finite time in $L^\infty$ ! In particular, classical bounded solutions do not exist globally in time.
\vskip5mm
If  $u_0,v_0, w_0\in L^\infty(\Omega)$, local existence and uniqueness of nonnegative and uniformly bounded solution to $(\mathcal S)$ are known (see e.g. \cite{R}). More precisely, there exists $T>0$ and a unique classical  solution $(u,v,w)$ of $(\mathcal{S})$ on $[0,T)$. If $T_{\max}$ denotes the greatest of these T's, then
\begin{equation}
\big(T_{\max}<+\infty\big)\Longrightarrow\lim_{t\nearrow T_{\max}}\left(\Vert u(t)\Vert_{L^\infty(\Omega)}+\Vert v(t)\Vert_{L^\infty(\Omega)}+ \Vert w(t)\Vert_{L^\infty(\Omega)}\right)=+\infty.
\end{equation}
To prove global existence ($i.e.\;T_{\max}=+\infty$), it is sufficient to obtain an a priori estimate of the form
\begin{equation}
\forall t\in[0,T_{\max}),\qquad\Vert u(t)\Vert_{L^\infty(\Omega)}+\Vert v(t)\Vert_{L^\infty(\Omega)}+\Vert w(t)\Vert_{L^\infty(\Omega)}\leq H(t),
\end{equation}
where  $H : [0,+\infty)\rightarrow [0,+\infty)  $ is a nondecreasing and continuous function.
\vskip5mm
This type of estimates is far of being obvious for our system
except the case where diffusion coefficients $d_1$, $d_2$, $d_3$ are equal $i.e\; d_1=d_2=d_3=d$. Indeed,  $Z=u+v+2w$ satisfies
$$(E)\left\lbrace\begin{array}{llll}
Z_t-d\Delta Z & = & 0& (0,+\infty)\times\Omega,\\
\Part{Z}{n}   & = & 0 & (0,+\infty)\times\partial\Omega,\\
Z(0,x) & = & Z_0(x) & x\in\Omega,\\
\end{array}
\right.$$
where $Z_0(x)=u_0(x)+v_0(x)+2w_0(x)$.\\
In particular, we deduce by maximum principle that
$$\|u(t)+v(t)+2w(t)\|_{L^\infty(\Omega)}\leq \|u_0+v_0+2w_0\|_{L^\infty(\Omega)},\quad\quad t\geq 0.$$
Together with nonnegativity, this implies
 $$\|u(t)\|_{L^\infty(\Omega)}+ \|v(t)\|_{L^\infty(\Omega)}+ \|w(t)\|_{L^\infty(\Omega)}\leq \|u_0+v_0+2w_0\|_{L^\infty(\Omega)}, \quad\quad t\geq 0.$$
In other words, $u(t)$, $v(t)$ and $w(t)$ stay uniformly bounded in $L^\infty(\Omega)$ and therefore $T_{\max}=+\infty$.
\vskip2mm
In the case where the diffusion coefficients are different from each other, global existence  is considerably  more complicated. It has been studied by several authors in the following cases.\\
{\bf First case } $\alpha=\beta=\gamma=1$.
 \par In this case, global existence of classical solutions has been obtained by Rothe \cite{R} for dimension $N\leq 5$. Later, it has first been proved by Pierre \cite{P1} for all dimensions $N$ and then by Morgan  \cite{M}.
 \par The exponentiel decay towards equilibrium has been studied by Desvillettes-Fellner \cite{DF} in the case of one space dimension.
 \par The global existence of weak solutions  has been proved by Laamri \cite{L} for initial data $u_0$,  $v_0$ and $w_0$ only in $L^1(\Omega)$.\\
\noindent{\bf Second case }$\gamma = 1$ regardless of $\alpha$ and $\beta$.
 \par In this case, global existence of classical solutions has been obtained by Feng \cite{F} in  all dimensions $N$ and  more general boundary conditions.\\
\noindent{\bf Third case } $\alpha+\beta\leq 2$ or $\gamma\leq 2$.
\par In this  case, Pierre \cite{P2} has proved global existence of weak solutions for initial data $u_0$,  $v_0$ and $w_0$ only in $L^2(\Omega)$.
\vskip3mm
Our paper mainly completes the investigations of [\cite{F}, \cite{M}, \cite{P1}, \cite{R}] and [\cite{L}, \cite{P2}]. As far as we know, our results are new either when $\alpha+\beta< \gamma$, or when $1<\gamma<\dis\frac{N+6}{N+2}$ regardless of $\alpha$ and $\beta$. For the sake of clarity, we decided to focus in this work on the question of global existence in time of solutions in the case of homogeneous Neumann boundary conditions. So, we shall prove global existence of classical solutions to system $(\mathcal S)$ in the following cases :\\
* $\alpha+\beta< \gamma$ ; \\
* ($d_1=d_3$ or $d_2=d_3$) and for any $(\alpha,\beta,\gamma)$ ;\\
* $d_1=d_2$ and for any  $(\alpha,\beta,\gamma)$ such that $\alpha+\beta\neq\gamma$ ; \\
* $1<\gamma < \dis\frac{N+6}{N+2}$ and for any $(\alpha,\beta)$.
\par For the sake of completeness and for the reader's convenience, we shall also give a direct proof different from that of Feng \cite{F} in the special case $\gamma = 1$.
\vskip3mm
\noindent{\bf Notation : } Throughout this study, we denote by $C_i$'s various positive numbers depending only on the data and  for $p\in [1,+\infty[$
$$\|u(t)\|_p=\left(\int_\Omega |u(t,x)|^p\,dx\right)^{1/p}, \quad \|u\|_{L^p(Q_T)}=\left(\int_0^T\int_\Omega |u(t,x)|^p\,dt dx\right)^{1/p},$$
$$\|u(t)\|_\infty=\dis\text{esse sup}_{x\in \Omega} |u(t,x)|, \quad\|u\|_{L^\infty(Q_T)}=\dis\text{esse sup}_{(t,x)\in Q_T} |u(t,x)|.$$
\section{The main results}
One of the main ingredients for the proof of our results is the following lemma which is based on the regularizing effects of the heat equation. This lemma has been introduced by  Hollis-Martin-Pierre in \cite{HMP}.
\begin{lem}
Let $T>0$ and $(\phi,\psi)$ the classical solution of
$$\left\lbrace\begin{array}{llll}
\phi_t-d_1\Delta \phi &=& f(\phi,\psi) & (t,x)\in (0,T)\times\Omega\\
\psi_{t}-d_2\Delta \psi &=& g(\phi,\psi) &  (t,x)\in (0,T)\times\Omega\\
\Part{\phi}{n}(t,x)& = &0 &(t,x)\in (0,T)\times\partial\Omega\\
\Part{\psi}{n}(t,x)& = &0 & (t,x)\in (0,T)\times\partial\Omega\\
\phi(0,x) & = & \phi_0(x) & x\in\Omega\\
\psi(0,x) & = & \psi_0(x) & x\in\Omega.
\end{array}
\right. $$
Assume that $f+g=0$, then for each $p\in(1,+\infty)$, there  exists  $C$ such that for all $t\in (0,T)$
\begin{equation} \label{inegalitefondamentale1}
\|\psi\|_{L^p(Q_t)}\leq C\left[\|\phi\|_{L^p(Q_t)}+1\right].
\end{equation}
\end{lem}
A more general version of this lemma can be founded in \cite[lemma 3.4]{P2}.   $\square$
\subsection{The case $\alpha+\beta<\gamma$}
%
%
\begin{thm} Assume that $0\leq u_0, v_0, w_0\leq M$ where $M$ is a positive real.\\
If  $\alpha+\beta<\gamma$, then the system $(\mathcal S)$ admits a global classical solution.
\end{thm}
{\bf Proof } :\\
$\bullet$  Let $T\in (0, T_{\max})$ and  let $t\in (0,T]$.
 Thanks to the nonnegativity of $u$, $v$ and $w$, we deduce from the equation $(E_1)$ that $u$ is bounded from above by the solution $U$ of
$$(P_1)\left\lbrace\begin{array}{llll}
U_t-d_1\Delta U & = & w^\gamma & (t,x)\in (0,T)\times\Omega\\
\Part{U}{n}(t,x)   & = & 0 &  (t,x)\in(0,T)\times\partial\Omega\\
U(0,x) & = & u_0(x) & x\in\Omega,
\end{array}
\right.$$
and we deduce from the equation $(E_2)$ that $v$ is bounded from above by the solution $V$ of
$$(P_2)\left\lbrace\begin{array}{llll}
V_t-d_2\Delta V & = & w^\gamma &(t,x)\in (0,T)\times\Omega\\
\Part{V}{n}(t,x)   & = & 0 & (t,x)\in(0,T)\times\partial\Omega\\
V(0,x) & = & v_0(x) & x\in\Omega.
\end{array}
\right.$$
Therefore it is sufficient to show that $ w \in L^p(Q_T)$ for $p$ large enough.\\
$\bullet$
Let $q>1$. Multiplying the equation $(E_3)$ by $w^q$ and integrating over $Q_T$, we get
\begin{equation}\label{integrationapresmultiplicationparwq} 
\frac{1}{q+1}\int_\Omega w^{q+1}(T) +qd_3\int\int_{Q_T}|\nabla w |^2w^{q-1}+\int\int_{Q_T}w^{q+\gamma}= \int\int_{Q_T} u^\alpha v^\beta w^q+ K_0
\end{equation}
where%
$$K_0= \frac{1}{q+1}\int_\Omega w^{q+1}_0.$$
Thanks to H\"{o}lder's inequality, we have
\begin{equation}\label{Holder1} 
\int\int_{Q_T} u^\alpha v^\beta w^q\leq \|u\|^\alpha_{L^{\alpha r}(Q_T)} \|v\|^\beta_{L^{\beta s}(Q_T)} \|w\|^q_{L^{\gamma+ q}(Q_T)}
\end{equation}
where
$$\frac{1}{r}+\frac{1}{s}+\frac{q}{q+\gamma}=1.$$
Since $\alpha +\beta<\gamma$, we can choose   $r$ such that   $r\alpha \leq q+\gamma$ and $s$ such that $s\beta \leq q+\gamma$. To convince oneself, it is enough to draw the straight line with cartesian equation $x+y=\dis\frac{\gamma}{q+\gamma}$ and to identify the points with coordinates $(\dis\frac{\alpha}{q+\gamma},0)$ and  $(0, \dis\frac{\beta}{q+\gamma})$.\\
Then  $ L^{q+\gamma}(Q_T)\subset  L^{\alpha r}(Q_T)$ and $L^{q+\gamma}(Q_T)\subset L^{\beta s}(Q_T) $. Consequently, there exists
  $C_1$ 
 such that
\begin{equation}\label{LpinclusdansdansLqdesquepinferieurq} 
\int\int_{Q_T} u^\alpha v^\beta w^q\leq C_1\|u\|^\alpha_{L^{\gamma+ q}(Q_T)} \|v\|^\beta_{L^{\gamma+ q}(Q_T)} \|w\|^q_{L^{\gamma+ q}(Q_T)}.
\end{equation}
By virtue of  lemma 1,  there  exists $C_2$
such that
\begin{equation}\label{equivalencedesnormes1} 
\|u\|_{L^{\gamma+q}(Q_T)}\leq  C_2(1+\|w\|_{L^{\gamma+q}(Q_T)})
\end{equation}
and there exists  $C_3$ 
such that
\begin{equation}\label{equivalencedesnormes2}
\|v\|_{L^{\gamma+q}(Q_T)}\leq C_3(1+\|w\|_{L^{\gamma+q}(Q_T)}).
\end{equation}
Thanks to (\ref{equivalencedesnormes1}) and (\ref{equivalencedesnormes2}), estimate (\ref{LpinclusdansdansLqdesquepinferieurq}) can be written
\begin{equation} 
\int\int_{Q_T} u^\alpha v^\beta w^q\leq C_4\left(1+\|w\|_{L^{\gamma+ q}(Q_T)}\right)^{\alpha}\left(1+\|w\|_{L^{\gamma+ q}(Q_T)}\right)^{\beta}\left(1+\|w\|_{L^{\gamma+ q}(Q_T)}\right)^{q}.
\end{equation}
If $\|w\|_{L^{\gamma+ q}(Q_T)}\leq 1$ then the proof ends up. Otherwise, there exists $C_5$ such that
\begin{equation}
\int\int_{Q_T} u^\alpha v^\beta w^q\leq C_5\|w\|^{q+\alpha+\beta}_{L^{\gamma+ q}(Q_T)}.
\end{equation}
So we deduce from (\ref{integrationapresmultiplicationparwq})
\begin{equation}\label{21}
\int\int_{Q_T}w^{q+\gamma}\leq C_5\|w\|^{q+\alpha+\beta}_{L^{\gamma+ q}(Q_T)}+ K_0.
\end{equation}
%
With the notation $R:= \dis\int\int_{Q_T} w^{q+\gamma}$, estimate (\ref{21})  can be written
\begin{equation}\label{22}
R\leq C_5 R^{\frac{q+\alpha+\beta}{q+\gamma}}+K_0.
\end{equation}
Since  $q+\alpha+\beta < q+ \gamma$, by applying Young's inequality  to (\ref{22}), we obtain
\begin{equation}\label{23}
(1-\varepsilon)R\leq  K_0 +C_6.
\end{equation}
Then, for   $\varepsilon\in (0,1)$, we have the desired estimate
\begin{equation}\label{24}
\|w\|_{L^{q+\gamma}(Q_T)}\leq C_{7}.
\end{equation}
Going back to $(P_1)$ and $(P_2)$, we have,  by choosing $q$ such that $\dis\frac{q+\gamma}{\gamma}>\dis\frac{N+2}{2}$ and thanks to the $L^p$-regularity theory for the heat operator (see \cite{LSU}),
  \begin{eqnarray}\label{25}
\|u\|_{L^\infty(Q_T)}&\leq & C_8
\end{eqnarray}
\begin{eqnarray}\label{26}
\|v\|_{L^\infty(Q_T)}&\leq & C_{9}.
\end{eqnarray}
Now going  back to $(E_3)$, we deduce from (\ref{25}) and (\ref{26}) that there exists  $C_{10} $ such that
\begin{eqnarray}
\|w\|_{L^\infty(Q_T)}\leq  C_{10}.
\end{eqnarray}
This implies that $T_{\max}=+\infty$.   $\square$
\vskip3mm
\noindent{\bf Remark } This method seems to be  specific to the case $\alpha+\beta<\gamma$. It fails when $\alpha+\beta\geq \gamma$ since some restrictions on the parameters $\alpha$, $\beta$, $\gamma$ and on the diffusion coefficients  will appear.
\subsection{Case where $d_1=d_3$ or $d_2=d_3$ or  $d_1=d_2$.}
%
%
\begin{thm} Assume that $0\leq u_0, v_0, w_0\leq M$.\\
(i) If $d_1=d_3$ or $d_2=d_3$, then  system $(\mathcal S)$ admits a global classical solution for any  $(\alpha, \beta, \gamma)$.\\
(ii) If $d_1=d_2$, then the system $(\mathcal S)$  admits a global classical solution for any $(\alpha, \beta, \gamma)$ such that $\alpha+\beta\neq\gamma$.
\end{thm}
{\bf Proof} :\\
$(i)$ Assume that $d_1=d_3=d$, we have
$$(u+w)_t -d\Delta(u+w)=0\;;\;\Part{(u+w)}{n}=0\;;\;(u+w)(0,x)=u_0(x)+w_0(x) .$$
We deduce by maximum principle
\begin{equation}\label{28}
 \|u(t)+w(t)\|_\infty\leq \|u_0+w_0\|_\infty.
\end{equation}
Together with the nonnegativity of $u$ et $w$, this implies that $u(t)$ and $w(t)$ are uniformly bounded in $L^\infty(\Omega)$.\\
By going  back to $(E_2)$ and thanks to the $L^p$-regularity theory for the heat operator (see \cite{LSU}), we conclude
 that   $\|v(t)\|_\infty$ is uniformly bounded in $L^\infty(\Omega)$ on all interval $[0,T]$ so that $T_{\max}=+\infty$.\\
$(ii)$ Assume that $d_1=d_2=d$. The case  $\alpha+\beta<\gamma$ was already handled in the theorem 1, so it remains only to tackle the case $\gamma< \alpha+\beta$. Moreover, one can assume that $u_0\neq v_0$ since if  $u_0=v_0$ the result is obvious.\\
Since $d_1=d_2=d$, we have
$$(u-v)_t -d\Delta(u-v)=0\;;\;\Part{(u-v)}{n}=0\;;\;(u-v)(0,x)=u_0(x)-v_0(x) .$$
The maximum principle then implies  $\|u(t)-v(t)\|_\infty\leq \|u_0-v_0\|_\infty=C $.
 Hence we have
\begin{eqnarray*}
u^{\alpha+\beta}&=&u^{\alpha}v^\beta +u^{\alpha}(u^{\beta}-v^{\beta} ) \\
&=&u^{\alpha}v^\beta +u^{\alpha}\beta(\theta u + (1-\theta)v )^{\beta-1}(u-v)\text{ where } \theta\in ]0,1[\\
&\leq&u^{\alpha}v^\beta +u^{\alpha}\beta2^{\beta-1}C(u^{\beta-1}+v^{\beta-1}).
\end{eqnarray*}
Thanks to Young's inequality, there exists  $C_{11}>0$ and  $C_{12}>0$ such that
\begin{eqnarray}\label{29}
C_{11}u^{\alpha+\beta}\leq u^{\alpha}v^\beta +C_{12}.
\end{eqnarray}
By virtue of (\ref{29}),  equation $(E_1)$ implies that
\begin{eqnarray}\label{30}
u_t-d_1\Delta u+C_{11} u^{\alpha+\beta}\leq w^\gamma +C_{12}.
\end{eqnarray}
Let $q>1$. Multiplying (\ref{30}) by  $u^q$ and integrating over $Q_T$, we obtain
\begin{equation}\label{integrationapresmultiplicationparuq}
\frac{1}{q+1}\int_\Omega u^{q+1}(T) +qd_2\int\int_{Q_T}|\nabla u |^2u^{q-1}+C_{11}\int\int_{Q_T}u^{q+\alpha+\beta}\leq\int\int_{Q_T} w^\gamma u^q+ C_{12}\int\int_{Q_T}u^q + K_1
\end{equation}
where%
$$K_1= \frac{1}{q+1}\int_\Omega u^{q+1}_0.$$
Thanks to H\"{o}lder's inequality, we have
\begin{equation}\label{Holder2}
\int\int_{Q_T}  w^\gamma u^q \leq \left(\int\int_{Q_T} w^{\gamma r}\right)^{1/r}\left(\int\int_{Q_T} u^{qs}\right)^{1/s}
\end{equation}
where $r=\dis\frac{\alpha+\beta+ q}{\gamma}$ and $s=\dis\frac{\alpha+\beta+ q}{q+\alpha+\beta-\gamma}$.\\
Lemma 1 implies that there exists $C_{13}$ such that
\begin{eqnarray*}
\left(\int\int_{Q_T} w^{\gamma r}\right)^{1/r}= \|w\|_{L^{q+\alpha+\beta}(Q_T)}^\gamma\leq C_{13}^\gamma\left(1+\|u\|_{L^{q+\alpha+\beta}(Q_T)}\right)^\gamma.
\end{eqnarray*}
If $\|u\|_{L^{q+\alpha+ \beta}(Q_T)}\leq 1$ then the proof ends up. Otherwise, there exists $C_{14}$ such that
\begin{equation}\label{dominationdewparu}
\left(\int\int_{Q_T} w^{\gamma r}\right)^{1/r}\leq C_{14}\|u\|^\gamma_{L^{q+\alpha+\beta}(Q_T)}.
\end{equation}
Since  $qs<q+\alpha+\beta$, we have $L^{q+\alpha+\beta}(Q_T)\subset L^{qs}(Q_T)$, then there exists $C_{15}$ such that
\begin{equation}\label{LqsInclusLqplusalphaplusbeta}
\left(\int\int_{Q_T}u^{qs}\right)^{1/s}\leq C_{15}\|u\|^q_{L^{q+\alpha+\beta}(Q_T)}.
\end{equation}
Denote $S: =\dis\int\int_{Q_T}u^{q+\alpha+\beta}$. Estimates (\ref{dominationdewparu}) and (\ref{LqsInclusLqplusalphaplusbeta}) imply that
\begin{equation}
\int\int_{Q_T}w^\gamma u^q\leq C_{16}S^{\frac{q+\gamma}{q+\alpha+\beta}}.
\end{equation}
Moreover, since  $ L^{q+\alpha+\beta}(Q_T)\subset L^{q}(Q_T)$, there exists $C_{17}$ such that
\begin{equation}\label{LqInclusLqplusalphaplusbeta}
C_{12}\int\int_{Q_T}u^q\leq C_{17}S^{\frac{q}{q+\alpha+\beta}}.
\end{equation}
Since $\gamma <\alpha+\beta$, by applying Young's inequality, there exists $C_{18}$  such that
\begin{equation}\label{Young1}
C_{16}S^{\frac{q+\gamma}{q+\alpha+\beta}}\leq \frac{\ve}{2}S+C_{18}.
\end{equation}
Applying again Young's inequality, there exists $C_{19}$  such that
\begin{equation}\label{Young2}
C_{17}S^{\frac{q}{q+\alpha+\beta}}\leq \frac{\ve}{2}S+C_{19}.
\end{equation}
Consequently, estimate (\ref{integrationapresmultiplicationparuq}) implies
\begin{equation}\label{EstimateEpsilon}
(C_{11}-\ve)S\leq C_{18}+C_{19}+K_1.
\end{equation}
By choosing $\ve<C_{11}$ in (\ref{EstimateEpsilon}), there exists $C_{20}$ such that
\begin{equation}\label{FinalEstimateonu}
\|u\|_{L^{q+\alpha+\beta}(Q_T)}\leq C_{20}.
\end{equation}
Thanks to lemma 1 and estimate (\ref{FinalEstimateonu}) there exists $C_{21}$ such that
\begin{equation}\label{FinalEstimateonw}
\|w\|_{L^{q+\alpha+\beta}(Q_T)}\leq C_{21}.
\end{equation}
By going back to $(P_1)$ and $(P_2)$, we have  by choosing $q$ such that $\dis\frac{q+\alpha+\beta}{\gamma}>\dis\frac{N+2}{2}$ and thanks to the $L^p$-regularity theory for the heat operator (see \cite{LSU})
  \begin{eqnarray}\label{EstimateLinfinieOnu}
\|u\|_{L^\infty(Q_T)}&\leq & C_{22}
\end{eqnarray}
\begin{eqnarray}\label{EstimateLinfinieOnv}
\|v\|_{L^\infty(Q_T)}&\leq & C_{23}.
\end{eqnarray}
Now let's go  back to $(E_3)$, we deduce from (\ref{EstimateLinfinieOnu}) and (\ref{EstimateLinfinieOnv}) that there exists  $C_{24} $ such that
\begin{eqnarray}
\|w\|_{L^\infty(Q_T)}\leq  C_{24}.
\end{eqnarray}
This implies that $T_{\max}=+\infty$.   $\square$\\
{\bf Remark } Even in the last  case $i.e\; d_1=d_2$, global existence or blow-up in the limit case  $\alpha+\beta = \gamma$ remain an open problem. $\square$
\subsection{Case $1\leq \gamma< \dis\frac{N+6}{N+2}$ regardless of $\alpha$ and $\beta$.}
%
%
\begin{thm}
Assume that $0\leq u_0, v_0, w_0\leq M$ where $M>0$.
If $1\leq \gamma<\dis\frac{N+6}{N+2}$, then the system $(\mathcal S)$ admits a global classical solution for any  $(\alpha,\beta)\in [1,+\infty)^2$.
\end{thm}
{\bf Proof } :\\
 Let $T\in (0, T_{\max})$ and  let $t\in (0,T]$. Thanks to the nonnegativity of $u$, $v$ and $w$, we deduce from the equation $(E_1)$ that $u$ is bounded from above by the solution $U$ of
$$(P_1) \left\lbrace\begin{array}{llll}
{U}_{t}-d_1\Delta U &=& w^\gamma & (t,x)\in (0,T)\times\Omega\\
\Part{U}{n}(t,x)& = &0 & (t,x)\in (0,T)\times\partial\Omega\\
U(0,x) & = & u_0(x) & x\in\Omega.
\end{array}
\right. $$
Therefore it is sufficient to show that $ w \in L^p(Q_T)$ for $p$ large enough. For this we have to distinguish the case $\gamma=1$ and the case $\gamma>1$.\\
$\bullet$ {\bf Case $\gamma=1$ and $\alpha, \beta\geq 1$.}
%
\par  Let us recall that global existence of classical solutions for $(\mathcal S)$ when $\alpha=\beta=\gamma=1$ has been studied by several authors. It has been obtained by  Rothe \cite{R} for dimension $N\leq 5$. Later, it has first been proved by Pierre \cite{P1} for all dimensions $N$ and then by Morgan  \cite{M}.\\
 Independantly,  Feng \cite{F} has proved global existence  in the case $\gamma=1$ regardless of $\alpha$ and $\beta$ and more general boundary conditions.\\
For the sake of completeness and for the reader's convenience, we give here a simple and direct proof in the last case ($\gamma=1$ regardless of $\alpha$ and $\beta$). In our proof, we use an idea introduced by Pierre in \cite{P1} and applied in \cite{MP}.\\
\par For any $p\geq 1$, we deduce from $(P_1)$ and the semigroup property
\begin{equation}\label{ineq1}
\|u(t)\|_p\leq \|u_0\|_p+\int_0^t\|w(s)\|_p\,ds.
\end{equation}
By applying H\"{o}lder's inequality for $p> 1$ and thanks to (\ref{inegalitefondamentale1}), we obtain
\begin{equation}\label{ineq2}
\int_0^t\|w(s)\|_p\,ds\leq t^{1/p'}\left(\int_0^t\int_\Omega w^p\,dsdx\right)^{1/p}\leq t^{1/p'}C_{25}\left[1+\left(\int_0^t\int_\Omega u^p\,dsdx\right)^{1/p}\right]
\end{equation}
where $p'=\dis\frac{p}{p-1}$.\\
For $t\in (0,T]$, let us set $h(t) :=\displaystyle\int_\Omega |u(t,x)|^p\,dx$. Inequality (\ref{ineq1}) can be written
\begin{equation}\label{ineq3}
h(t)^{1/p}\leq C_{26}+ C_{27}\left(\int_0^t h(s)\,ds\right)^{1/p}.
\end{equation}
%
Taking the $p^{\text{th}}$ power of  (\ref{ineq3}) 
we obtain
%
\begin{equation} \label{ineq4}
h(t)\leq 2^{p-1}C_{26}^p +2^{p-1}C_{27}^p\int_0^t h(s)\,ds.
\end{equation}
But, inequality  (\ref{ineq4}) is a linear Gronwall's inequality, then
\begin{equation}\label{inegalite7}
\|u\|_{L^p(Q_T)}\leq C_{28}.
\end{equation}
Repeating the method above with $v$ instead of $u$, we obtain
\begin{equation}\label{inegalite8}
\|v\|_{L^p(Q_T)}\leq C_{29}.
\end{equation}
Estimates (\ref{inegalite7}) and (\ref{inegalite8}) imply that for some $q>\dis\frac{N+2}{2}$
\begin{equation}\label{inegalite9}
\|u^\alpha v^\beta\|_{L^q(Q_T)}\leq C_{30}.
\end{equation}
Going back to equation  $(E_3)$ we have, thanks to the $L^q$-regularity theory for the heat operator,
\begin{equation}\label{inegalite10}
\|w\|_{L^\infty(Q_T)}\leq  C_{31}.
\end{equation}
This concludes the proof for the case $\gamma=1$ regardless of $\alpha$ and $\beta$.   $\square$
\vskip2mm
\noindent$\bullet$ {\bf Case $1<\gamma<\dis\frac{N+6}{N+2}$. }\\
The proof in this case is based on lemma 1 and these two following lemmas.
\begin{lem}[Michel Pierre]\label{EstimationLdeux}
 Let $T>0$ and let $Z$ the  solution of
$$\left\lbrace\begin{array}{llll}
Z_t-\Delta(A(t,x) Z) & \leq  & 0 & (t,x)\in (0,T)\times\Omega,\\
\Part{Z}{n}(t,x)   & = & 0 & (t,x)\in (0,T)\times\partial\Omega,\\
Z(0,x) & = & Z_0(x) & x\in\Omega.
\end{array}
\right.$$
Assume that $0<d<A(t,x)<D$ where $(d,D)\in (0,+\infty)^2$. Then, there exists $C=C(T,d,D,\Omega)$ such that
 $$\|Z\|_{L^2(Q_T)}\leq C\|Z_0\|_{L^2(\Omega)}.$$
\end{lem}
\par For a general version of this  lemma, see  \cite[proposition 6.1]{P2} or \cite[theorem 3.1]{DFPV}. $\square$
\vskip2mm
\begin{lem}\label{estimationLpLq}
Let $(p,q)$ such that $1\leq p\leq q\leq +\infty$, $d>0$ and $S_d(t)$ the semigroup generated in $L^p(\Omega)$  by $-d\Delta$ with homogeneous Neumann boundary condition.
Then
\begin{equation}\label{EstimationLpLq}
\|S_d(t)Y\|_q\leq \left(C(\Omega)m(t)\right)^{\frac{-N}{2}(\frac{1}{p}-\frac{1}{q})}\|Y\|_p, \text{ for all } Y \in L^p(\Omega),\;\; t>0
\end{equation}
where $m(t)=\min(1,t)$.
\end{lem}
For a proof of this lemma see  for instance  \cite[Lemma 3, p. 25]{R} or \cite[Theorem 3.2.9, p. 90]{D}.   $\square$
\vskip3mm
We now go back to the proof of theorem 3.\\

By applying lemma 2 to the system $(\mathcal{S})$ where $Z=u+v+2w$ and $A=\dis\frac{d_1u+d_2v+2d_3w}{u+v+2w}$, we have $ u,\; v,\; w\; \in L^2(Q_T)$. More precisely, there exists $C_{32}$ such that
\begin{equation}\label{estimationL2Ofuandvandw}
\|u\|_{L²(Q_T)},\,\|v\|_{L²(Q_T)},\,\|w\|_{L²(Q_T)} \leq C_{32}.
\end{equation}

Now, we have thanks to the estimate (\ref{EstimationLpLq}) with $p>1$ and  $q=+\infty$
\begin{eqnarray}\label{Thm3.1}
\|u(t)\|_\infty &\leq& \|u_0\|_\infty+ C_{33}\int_0^t (t-s)^{\frac{-N}{2p}}\|w^\gamma(s) \|_p\,ds.
\end{eqnarray}
By applying H\"{o}lder's inequality, we obtain
\begin{eqnarray}\label{Thm3.2}
\int_0^t (t-s)^{\frac{-N}{2p}}\|w(s)^\gamma \|_p\,ds\leq \left(\int_0^t (t-s)^{\frac{-Np'}{2p}}\,ds\right)^{1/p'}
\left(\int_0^t \|w^\gamma(s)\|_p^p\,ds\right)^{1/p}.
\end{eqnarray}
We first remark that the integral $\dis\int_0^t(t-s)^{\frac{-N}{2(p-1)}}ds$ converges
when $p>\dis\frac{N+2}{2}$ and we have
\begin{eqnarray*}
\int_0^t(t-s)^{\frac{-Np'}{2p}}\,ds &=& t^{1-N/(2(p-1))}\int_0^1(1-y)^{\frac{-N}{2(p-1)}}\,dy\\
&\leq &C(T)^{p/(p-1)}= T^{1-N/(2(p-1))}\int_0^1(1-y)^{\frac{-N}{2(p-1)}}\,dy.
\end{eqnarray*}
On the other hand, lemma 1 implies that
\begin{eqnarray}\label{Thm3.3}
\left(\int_0^t \|w^\gamma(s)\|_p^pds\right)^{1/p}=\|w\|_{L^{p\gamma}(Q_t)}^\gamma\leq C_{34}^\gamma\left(1+\|u\|_{L^{p\gamma}(Q_t)}\right)^\gamma.
\end{eqnarray}
If $\|u\|_{L^{p\gamma}(Q_t)}\leq 1$ then the proof ends up. Otherwise there exists $C_{35}$ such that
\begin{eqnarray}\label{Thm3.4}
\left(\int_0^t \|w^\gamma(s)\|_p^pds\right)^{1/p}\leq C_{35}\|u\|_{L^{p\gamma}(Q_T)}^\gamma.
\end{eqnarray}
Since
\begin{eqnarray*}\label{Thm3.5}
\|u\|_{L^{p\gamma}(Q_T)}^\gamma = \left(\int\int_{Q_T}u^{p\gamma}\right)^{1/p}=
\left(\int\int_{Q_T}u^{p\gamma-p+\ve + p-\ve}\right)^{1/p}\\
\leq\|u\|_{L^\infty{(Q_T)}}^{1-\ve/p}\left(\int\int_{Q_T}u^{p\gamma-p+\ve}\right)^{1/p},
\end{eqnarray*}
it follows that (\ref{Thm3.1}) can be written
\begin{eqnarray}\label{Thm3.6}
\|u(t)\|_\infty &\leq& \|u_0\|_\infty+ C_{36}\|u\|_{L^\infty{(Q_T)}}^{1-\ve/p}\left(\int\int_{Q_T}u^{p\gamma-p+\ve}\right)^{1/p}.
\end{eqnarray}
If $p(\gamma-1)<2$, by  choosing $\ve \in (0,\min(p, 2-p(\gamma-1))$, we deduce from (\ref{estimationL2Ofuandvandw}) and (\ref{Thm3.6}) that there exists $ C_{37}$ such that
\begin{equation}\label{UniformEstimateOfu}
\|u\|_{L^\infty(Q_T)}\leq C_{37}.
\end{equation}
Note that the above condition $p(\gamma-1)<2$ holds if $\gamma< 1+\dis\frac{2}{p}<1+\frac{4}{N+2}= \frac{N+6}{N+2}$.\\
We establish in the same  way that there exists  $C_{38}$ such that
\begin{equation}\label{UniformEstimateOfv}
\|v\|_{L^\infty(Q_T)}\leq C_{38}.
\end{equation}
Finally, for $(E_3)$, we deduce from (\ref{UniformEstimateOfu}) and (\ref{UniformEstimateOfv}) that there exists  $C_{39}$ such that
\begin{equation}\label{UniformEstimateOfw}
\|w\|_{L^\infty(Q_T)}\leq C_{39}.
\end{equation}
This concludes the proof in the case $1<\gamma<\dis\frac{N+6}{N+2}$.   $\square$\\
{\bf Remark } : Our conjecture is that  $\gamma^*=\dis\frac{N+6}{N+2}$ is not optimal. In fact, when $N=1$ one can prove that the result of theorem 3 still holds for $\gamma^*=7/2$.  $\square$
\section{Conclusion}
\noindent $\bullet$ All our results are still true if we replace homogeneous Neumann boundary conditions by homogeneous Dirichlet  boundary conditions, it suffices to replace  lemma 3 by the following one.
\begin{lem}
Let $(p,q)$ such that $1\leq p\leq q\leq +\infty$, $d>0$ and $S_d(t)$ the semigroup generated in $L^p(\Omega)$  by $-d\Delta$ with homogeneous Dirichlet boundary. Then
\begin{equation}\label{EstimationLpLqDirichlet}
\|S_d(t)Y\|_q\leq \left(4\pi t\right)^{\frac{-N}{2}(\frac{1}{p}-\frac{1}{q})}\|Y\|_p, \text{ for all } Y \in L^p(\Omega),\;\; t>0.
\end{equation}
\end{lem}
For a proof of this lemma, see  for instance  \cite[Proposition 48.4, p. 441]{QS}.   $\square$
\vskip3mm
\noindent $\bullet$ In the case where the diffusion coefficients are not equal (i.e. $d_i\neq d_j$ for all $1\leq i\neq j\leq 3$), global existence of classical solutions for $(\cal S)$ or  blow-up is still an open question when $$\dis\frac{N+6}{N+2} \leq \gamma \leq \alpha+\beta.$$    Our guess is that system $(\cal S)$ admits a classical global solution  for all  $\dis\frac{N+6}{N+2} \leq\gamma<\alpha+\beta$ and that there is a finite time blow-up when  $\gamma=\alpha+\beta$ and the  dimension $N$ is large.   $\blacksquare$
\vskip1cm
\noindent{\bf Acknowledgments :} I am indebted to Michel Pierre for numerous and helpful discussions concerning this work. I would like to thank warmly Philippe Souplet for his advice which improved the second part of the proof of theorem 3. I am sincerely grateful to Laurent Desvillettes for his careful reading of a initial version of this paper and for his useful comments. I thank Didier Schmitt  for many conversations.  Finally, I also thank the referee for valuable comments and remarks which allow me to improve the writing of this paper.

\end{document}